# INTERPLAY BETWEEN DIVIDEND RATE AND BUSINESS CONSTRAINTS FOR A FINANCIAL CORPORATION

By Tahir Choulli[1], Michael Taksar[2] and Xun Yu Zhou[3]

*University of Alberta, University of Missouri
and Chinese University of Hong Kong*

We study a model of a corporation which has the possibility to choose various production/business policies with different expected profits and risks. In the model there are restrictions on the dividend distribution rates as well as restrictions on the risk the company can undertake. The objective is to maximize the expected present value of the total dividend distributions. We outline the corresponding Hamilton–Jacobi–Bellman equation, compute explicitly the optimal return function and determine the optimal policy. As a consequence of these results, the way the dividend rate and business constraints affect the optimal policy is revealed. In particular, we show that under certain relationships between the constraints and the exogenous parameters of the random processes that govern the returns, some business activities might be redundant, that is, under the optimal policy they will never be used in any scenario.

**1. Introduction.** In recent years we have seen a lot of new results in the application of diffusion optimization models to financial mathematics. Together with portfolio optimization models, dividend distribution and risk control models have undergone major development.

In typical models of this type (see [[2], [3], [8]–[11], [13], [14], [16]]), the liquid assets of the company are governed by a Brownian motion with constant drift and diffusion coefficients. The drift term corresponds to the expected (potential) profit per unit time, while the diffusion term is interpreted as risk. The decrease of risk from business activities corresponds to a decrease

Received December 2002; revised March 2003.

[1]Supported by the Department of Systems Engineering and Engineering Management at the Chinese University of Hong Kong and NSERC Grant 121210818.

[2]Supported by NSF Grant DMS-00-72388.

[3]Supported by RGC Earmarked Grant CUHK 4054/98E.

*AMS 2000 subject classifications.* 91B70, 93E20.

*Key words and phrases.* Diffusion model, dividend distribution, business constraints, risk control, optimal stochastic control, HJB equation.

in potential profits. Different business activities in these models correspond to changing the drift and the diffusion coefficients of the underlying process *simultaneously*. This sets the scene for an optimal stochastic control model where the controls affect not only the drift, but also the diffusion part of the dynamic of the system.

In this article we study a model with an explicit restriction on risk control and on the rate at which the dividends are paid out. In addition, the company may have liability which it has to pay out at a constant rate no matter what the business plan is.

The controls are described by two functionals $a_t$ and $c_t$. The first represents the degree of business activity which the company assumes. The process $a_t$ takes on values in the interval $[\alpha, \beta]$, $0 < \alpha < \beta \leq +\infty$. The risk, which in our model is associated with the diffusion coefficient, and the potential profit, which is associated with the drift coefficient of the corresponding process, are both proportional to $a_t$. The constraints on the values of $a_t$ reflect institutional or statutory restrictions (e.g., for a public company) that the risk it can assume cannot exceed a certain level or that its business activities cannot be reduced to zero unless the company goes bankrupt.

The value $c_t$ of the second control functional shows the rate at which dividends are paid out at time $t$. The dividends are paid out from the liquid reserve and are distributed to shareholders. This corresponds to $c_t$ entering the drift coefficient of the reserve process with a negative sign. The dividend rate is bounded by a constant $M$ given a priori.

In our model we also assume the existence of a constant rate liability payment, such as a mortgage payment on a property or amortization of bonds. The results of this model can be viewed as an extension of the results of Choulli, Taksar and Zhou [4]. The presence of dividend rate constraints, however, adds a whole new dimension to the analysis as well as to the qualitative structure of the results obtained.

What is the most interesting is the interplay between the constraints and the exogenous parameters that govern the process of returns. Depending on the relationship between these parameters, we get several distinct cases of qualitative behavior of the company under the optimal policy.

This article is structured as follows. In the next section we present a rigorous mathematical formulation of the problem and state general properties of the *optimal return* or the *value* function. We also write the Hamilton–Jacobi–Bellman (HJB) equation this function must satisfy. In Section 3 we find a bounded smooth solution to the HJB equation. In Section 4 we construct the optimal policy and present our main findings in table form. Finally, in Section 5 we describe some economic interpretation of the results and state conclusions.

**2. Mathematical model.** We start with a filtered probability space $(\Omega, \mathcal{F}, \mathcal{F}_t, P)$ and a one-dimensional standard Brownian motion $W_t$ (with $W_0 = 0$) on it, adapted to the filtration $\mathcal{F}_t$. We denote by $R_t^\pi$ the reserve of the company at time $t$ under a control policy $\pi = (a_t^\pi, c_t^\pi; t \geq 0)$ (to be specified below). The dynamic of the reserve process $R_t^\pi$ is described by

$$(2.1) \qquad dR_t^\pi = (a_t^\pi \mu - \delta)\, dt + a_t^\pi \sigma\, dW_t - c_t^\pi\, dt, \qquad R_0^\pi = x,$$

where $\mu$ is the expected profit per unit time (profit rate), $\sigma$ is the volatility rate of the reserve process (in the absence of any risk control), $\delta$ represents the amount of money the company has to pay per unit time (the debt rate) irrespective of what business activities it chooses and $x$ is the initial reserve.

The control in this model is described by a pair of $\mathcal{F}_t$-adapted processes $\pi = (a_t^\pi, c_t^\pi; t \geq 0)$. A control $\pi = (a_t^\pi, c_t^\pi; t \geq 0)$ is admissible if $\alpha \leq a_t^\pi \leq \beta$ and $0 \leq c_t^\pi \leq M\ \forall t \geq 0$, where $0 < \alpha < \beta < +\infty$ and $0 < M < +\infty$ are given scalars. We denote the set of all admissible controls by $\mathcal{A}$. The control component $a_t^\pi$ represents one of the possible business activities available for the company at time $t$, and the component $c_t^\pi$ corresponds to the dividend payout rate at time $t$.

Given a control policy $\pi$, the time of bankruptcy is defined as

$$(2.2) \qquad \tau^\pi = \inf\{t \geq 0 : R_t^\pi = 0\}.$$

The *performance functional* associated with each control $\pi$ is

$$(2.3) \qquad J_x(\pi) = E\left( \int_0^{\tau^\pi} e^{-\gamma t} c_t^\pi\, dt \right),$$

where $\gamma > 0$ is an a priori given discount factor (used to calculate the present value of the future dividends), and the subscript $x$ denotes the initial state $x$. The objective is to find

$$(2.4) \qquad v(x) = \sup_{\pi \in \mathcal{A}} J_x(\pi)$$

and the optimal policy $\pi^*$ such that

$$(2.5) \qquad J_x(\pi^*) = v(x).$$

The exogenous parameters of the problem are $\mu$, $\sigma$, $\delta$, $\alpha$, $\beta$ and $\gamma$. The aim of this article is to obtain the optimal return function $v$ and the optimal policy *explicitly* in terms of these parameters.

The main tools for solving the problem are the dynamic programming and HJB equation (see [[6], [7] and [17] as well as relevant discussions in [[2], [9] and [16]]). We start by stating the following properties of the optimal return function $v$.

PROPOSITION 2.1. *The optimal return function $v$ is a concave, nondecreasing function subject to $v(0) = 0$ and*

$$0 \le v(x) \le \frac{M}{\gamma} \qquad \forall\, x > 0. \tag{2.6}$$

PROOF. The proof of the concavity and the monotonicity as well as the boundary condition, $v(0) = 0$, is similar to the one in [4]. To show (2.6), consider

$$0 \le E\left(\int_0^{\tau^\pi} e^{-\gamma t} c_t^\pi \, dt\right) \le M \int_0^\infty e^{-\gamma t} \, dt = \frac{M}{\gamma}. \qquad \square$$

If the optimal return function $v$ is twice continuously differentiable, then it must be a solution to the HJB equation

$$\begin{aligned}
0 &= \max_{\alpha \le a \le \beta, 0 \le c \le M} (\tfrac{1}{2}\sigma^2 a^2 V''(x) + (a\mu - \delta - c)V'(x) - \gamma V(x) + c) \\
&\equiv \max_{\alpha \le a \le \beta} (\tfrac{1}{2}\sigma^2 a^2 V''(x) + (a\mu - \delta)V'(x) - \gamma V(x) + M(1 - V'(x))^+)
\end{aligned} \tag{2.7}$$

$$V(0) = 0,$$

where $x^+ = \max(x, 0)$. This equation is rather standard and its derivation can be found in [[6], [7], and [17]]; see also [9] and [10].

Note that we do not know a priori whether the HJB equation has any solution other than the optimal return function. However, the following verification theorem, which says that *any* concave solution $V$ to the HJB equation (2.7) whose derivative is finite at 0 majorizes the performance functional for any policy $\pi$, is sufficient for us to identify optimal policies.

THEOREM 2.2. *Let $V$ be a concave, twice continuously differentiable solution of* (2.7), *such that $V'(0) < +\infty$. Then, for any policy $\pi = (a_t^\pi, c_t^\pi; t \ge 0)$,*

$$V(x) \ge J_x(\pi). \tag{2.8}$$

PROOF. Let $R_t^\pi$ be the reserve process given by (2.1). Denote the operator

$$L^a = \frac{1}{2}\sigma^2 a^2 \frac{d^2}{dx^2} + (a\mu - \delta)\frac{d}{dx} - \gamma.$$

Then by applying Ito's formula (see [5], Theorem VIII.27) to the process $e^{-\gamma t}V(R_t^\pi)$, we get

$$\begin{aligned}
e^{-\gamma(t\wedge\tau)}V(R_{t\wedge\tau}^\pi) &= V(x) + \int_0^{t\wedge\tau} e^{-\gamma s}\sigma a_s^\pi V'(R_s^\pi)\, dW_s \\
&\quad + \int_0^{t\wedge\tau} e^{-\gamma s} L^{a_s^\pi} V(R_s^\pi)\, ds - \int_0^{t\wedge\tau} e^{-\gamma s} V'(R_s^\pi) c_s^\pi \, ds.
\end{aligned} \tag{2.9}$$

Since $V$ is nondecreasing, concave with finite derivative at the origin, $V'(x)$ is bounded and the stochastic integral in (2.9) is a square integrable martingale whose expectation vanishes. In view of the HJB equation (2.7) and the inequality $c_s^\pi \le M$, we have

$$(2.10) \qquad L^{a_s^\pi} V(R_s^\pi) \le -c_s^\pi (1 - V'(R_s^\pi))^+.$$

Taking expectations of both sides of (2.9), in view of (2.10) we get

$$(2.11) \qquad \begin{aligned} &E(e^{-\gamma(t \wedge \tau)} V(R_{t \wedge \tau}^\pi)) \\ &\le V(x) - E \int_0^{t \wedge \tau} e^{-\gamma s} c_s^\pi [V'(R_s^\pi) + (1 - V'(R_s^\pi))^+] \, ds. \end{aligned}$$

Combining (2.11) with the fact that $y + (1 - y)^+ \ge 1$, we get

$$(2.12) \qquad E(e^{-\gamma(t \wedge \tau)} V(R_{t \wedge \tau}^\pi)) + E \int_0^{t \wedge \tau} e^{-\gamma s} c_s^\pi \, ds \le V(x).$$

Note that in view of the boundedness of $V'$,

$$e^{-\gamma(t \wedge \tau)} V(R_{t \wedge \tau}^\pi) \le e^{-\gamma t} K(1 + R_{t \wedge \tau}^\pi) \le e^{-\gamma t} K(1 + |R_t^\pi|)$$

for some constant $K$. Since $R_t^\pi$ is a diffusion process with uniformly bounded drift and diffusion coefficient, standard arguments yield $E|R_t^\pi| \le x + K_1 t$ for some constant $K_1$. Therefore,

$$(2.13) \qquad E e^{-\gamma(t \wedge \tau)} V(R_{t \wedge \tau}^\pi) \to 0$$

as $t \to \infty$. Thus taking the limit in (2.12) as $t \to \infty$, we arrive at

$$V(x) \ge E \int_0^\tau e^{-\gamma s} c_s^\pi \, ds = J_\pi(x). \qquad \square$$

The idea of solving the original optimization problem is first to find a concave, smooth function to the HJB equation (2.7) and then to construct a control policy [by solving a stochastic differential equation (SDE); for details see Section 4] whose performance functional can be shown to coincide with the solution to (2.7). Then, the above verification theorem establishes the optimality of the constructed control policy. As a by-product, there is no other concave solution to (2.7) than the optimal return function.

**3. A smooth solution to the HJB equation.** In this section, we are looking for a concave, smooth solution to (2.7). Assume that such a solution $V$ has been found. Let

$$(3.1) \qquad x_1 = \inf\{x \ge 0 : V'(x) \le 1\}.$$

Then, for $0 \le x < x_1$, (2.7) becomes

$$(3.2) \qquad 0 = \max_{\alpha \le a \le \beta} (\tfrac{1}{2} \sigma^2 a^2 V''(x) + (a\mu - \delta) V'(x) - \gamma V(x)),$$

while for $x \geq x_1$, (2.7) can be rewritten as

$$(3.3) \quad 0 = \max_{\alpha \leq a \leq \beta} (\tfrac{1}{2}\sigma^2 a^2 V''(x) + (a\mu - \delta - M)V'(x) - \gamma V(x) + M).$$

We start by seeking a smooth solution to (3.3). Obviously if $V'(0) \leq 1$, then $x_1 = 0$ and (2.7) is equivalent to (3.3) for all $x \geq 0$.

PROPOSITION 3.1.    If $\beta\mu \leq \delta$, then $V'(0) < 1$.

PROOF.    It follows from (2.7) that there exist $\tilde{a} \in [\alpha, \beta]$ such that

$$(3.4) \quad 0 = \tfrac{1}{2}\sigma^2 \tilde{a}^2 V''(0) + (\tilde{a}\mu - \delta)V'(0) + M(1 - V'(0))^+.$$

If $\beta\mu < \delta$, then each of the first two terms on the right-hand side of (3.4) is nonpositive with the second being strictly negative. Therefore, $M(1 - V'(0))^+ > 0$, which implies $V'(0) < 1$. The same argument goes if $\beta\mu = \delta$ and $V''(0) < 0$. In this case, either the first or the second term on the right-hand side of (3.4) is strictly negative. If $\beta\mu = \delta$ and $V''(0) = 0$, then the maximizer of the right-hand side of (2.7) is equal to $\beta$ for all $x$ in a right neighborhood of 0 [recall that $V'(0) > 0$]. Substituting $a = \beta$ either into (3.2) or into (3.3) and solving the resulting linear ordinary differential equation (ODE) with constant coefficients, we get a function $V$ whose second derivative at 0 does not vanish, which is a contradiction.   □

REMARK 3.2.    When the dividend rates are unrestricted, the condition $\beta\mu \leq \delta$ makes the problem trivial (see [4], Theorem 4.1). This is not the case when the dividend rates are bounded. Even if $\beta\mu \leq \delta$, the second derivative of $V$ at 0 is strictly negative, which makes the problem nontrivial in contrast to a similar situation in the case of unrestricted dividends.

Now we analyze the solution to (3.3) under the condition $\beta\mu > \delta$. As we see later, the qualitative nature of this solution depends on whether $a(x_1) < \alpha$ or $\alpha \leq a(x_1) < \beta$ or $a(x_1) \geq \beta$, where

$$(3.5) \qquad\qquad a(x) \equiv -\frac{\mu V'(x)}{\sigma^2 V''(x)} > 0, \qquad x < x_1.$$

To this end, we need the following proposition:

PROPOSITION 3.3.    (i) If $a(x_1) \geq \alpha$, then, for each $x \geq x_1$,

$$(3.6) \qquad\qquad\qquad a(x) \geq \alpha.$$

(ii) If $a(x_1) \geq \beta$, then, for each $x \geq x_1$,

$$(3.7) \qquad\qquad a(x) \geq \beta.$$

PROOF. (i) Suppose there exists $x_0 > x_1$ such that $a(x_0) < \alpha$. Then there exists $\varepsilon > 0$ such that $a(x) < \alpha$ for each $x$ with $|x-x_0| < \varepsilon$. Let $x' = \sup\{x_1 \leq x < x_0 : a(x) = \alpha\}$. Then $x_1 \leq x' < x_0 < x_0 + \varepsilon$ and $a(x') = \alpha$. Since $a(x) \leq \alpha$ for all $x \in [x', x_0 + \varepsilon)$, the function $V$ satisfies (3.3) with the maximum there attained at $a = \alpha$. Therefore,

$$(3.8) \quad V(x) = \frac{M}{\gamma} + K_1 \exp(\tilde{r}_+(\alpha)(x - x')) + K_2 \exp(\tilde{r}_-(\alpha)(x - x'))$$
$$\forall x \in [x', x_0 + \varepsilon).$$

Here

$$(3.9) \quad \tilde{r}_+(z) \equiv \frac{-(z\mu - \delta - M) + \sqrt{(z\mu - \delta - M)^2 + 2\gamma\sigma^2 z^2}}{z^2\sigma^2},$$

$$(3.10) \quad \tilde{r}_-(z) \equiv \frac{-(z\mu - \delta - M) - \sqrt{(z\mu - \delta - M)^2 + 2\gamma\sigma^2 z^2}}{z^2\sigma^2}, \qquad z > 0.$$

From (3.8) and (3.5), the equation $a(x') = \alpha$ can be rewritten as

$$K_1 \tilde{r}_+(\alpha) = -K_2 \tilde{r}_-(\alpha) \frac{\mu + \alpha\sigma^2 \tilde{r}_-(\alpha)}{\mu + \alpha\sigma^2 \tilde{r}_+(\alpha)},$$

which establishes a relationship between the constants $K_1$ and $K_2$. Using this relationship, we calculate

$$a(x) = \frac{-\mu V'(x)}{\sigma^2 V''(x)}$$

$$= \left( -\mu \left( \exp((\tilde{r}_+(\alpha) - \tilde{r}_-(\alpha))(x - x')) - \frac{\mu + \alpha\sigma^2 \tilde{r}_+(\alpha)}{\mu + \alpha\sigma^2 \tilde{r}_-(\alpha)} \right) \right)$$

$$(3.11) \qquad \times \left( \sigma^2 \left( \tilde{r}_+(\alpha) \exp((\tilde{r}_+(\alpha) - \tilde{r}_-(\alpha))(x - x')) \right. \right.$$

$$\left. \left. - \tilde{r}_-(\alpha) \frac{\mu + \alpha\sigma^2 \tilde{r}_+(\alpha)}{\mu + \alpha\sigma^2 \tilde{r}_-(\alpha)} \right) \right)^{-1}$$

$$\forall x \in [x', x_0 + \varepsilon).$$

However, we have $a(x) < \alpha$ for $x > x'$, which after a simple algebraic transformation of (3.11) is equivalent to $\exp((\tilde{r}_+(\alpha) - \tilde{r}_-(\alpha))(x - x')) < 1$. This leads to a contradiction. Therefore (3.6) holds.

(ii) By virtue of the assertion (i), $a(x) \geq \alpha$ for all $x \geq x_1$. Suppose there exists $x' > x_1$ such that $a(x') < \beta$. Then there exists $\varepsilon > 0$ such that $a(x) < \beta$ for all $x < x' + \varepsilon$. Let $\bar{x} = \sup\{x_1 \leq x < x' : a(x) = \beta\}$. Then $x_1 \leq \bar{x} < x'$ and $a(\bar{x}) = \beta$. In addition $\alpha \leq a(x) < \beta$ for all $\bar{x} < x \leq x'$. Substituting $a \equiv a(x)$ and $V''(x) = -(\mu V'(x))/(\sigma^2 a(x))$ into (3.3), we get

$$(3.12) \qquad 0 = \frac{\mu a(x)}{2} V'(x) - (\delta + M)V'(x) - \gamma V(x) + M.$$

Differentiating (3.12) and again substituting $V''(x) = -(\mu V'(x))/(\sigma^2 a(x))$ into the resulting equation, we obtain

$$(3.13) \qquad a(x)a'(x) = (a(x) - \tilde{c})\frac{\mu^2 + 2\sigma^2\gamma}{\mu\sigma^2}.$$

Then integrating (3.13) we get

$$(3.14) \quad a(x) - a(\bar{x}) + \tilde{c}\log\left(\frac{a(x) - \tilde{c}}{a(\bar{x}) - \tilde{c}}\right) = (x - \bar{x})\frac{\mu^2 + 2\sigma^2\gamma}{\mu\sigma^2} > 0$$
$$\forall \bar{x} < x < x',$$

which is a contradiction. Hence (3.7) holds and this completes the proof of the proposition. $\square$

First suppose that $a(x_1) \geq \beta$. In view of Proposition 3.3(i), we deduce that $a(x) \geq \beta$ for each $x \geq x_1$. Substituting $a = \beta$ into (3.3) and solving the resulting equation, we get

$$V(x) = \frac{M}{\gamma} + K_\beta \exp(\tilde{r}_-(\beta)(x - x_1)) \qquad \forall x \geq x_1.$$

Here $K_\beta$ is a constant which takes on either the value $-\frac{M}{\gamma}$ or $1/(\tilde{r}_-(\beta))$ depending, respectively, on whether $x_1$ in (3.1) is zero or not. Then straightforward calculations show that $a(x) = -\mu/(\sigma^2\tilde{r}_-(\beta))$. Thus, the condition $a(x_1) \geq \beta$ is equivalent to

$$(3.15) \qquad \tilde{c} \equiv \frac{2\mu(\delta + M)}{\mu^2 + 2\gamma\sigma^2} \geq \beta.$$

Next suppose $\alpha \leq a(x_1) < \beta$. By virtue of Proposition 3.3(i), $a(x) \geq \alpha$ for all $x \geq x_1$. As a result, $\alpha \leq a(x) < \beta$ in a right neighborhood of $x_1$. Substituting $a \equiv a(x)$ and $V''(x) = -(\mu V'(x))/(\sigma^2 a(x))$ into (3.3), we deduce that $a(x)$ satisfies (3.12). Then, following the same analysis there, we derive equation (3.13) for $a(x)$.

Suppose there exists $x' \geq x_1$ such that $a(x') < \tilde{c}$ [resp., $a(x') > \tilde{c}$]. Then from (3.13) we deduce that $a(x) < \tilde{c}$ [resp., $a(x) > \tilde{c}$] for each $x \geq x'$. Thus, by integrating (3.13) we derive (3.14) for all $x \geq x'$, with $\bar{x}$ replaced by $x'$. From (3.12) and (2.6), we see that $a(x) \leq \frac{2(\hat{\delta}+M)}{\mu}$ $\forall x \geq x_1$. Therefore, the left-hand side of (3.14) is bounded. This is a contradiction and we conclude that $a(x) = \tilde{c}$ for each $x \geq x_1$. In view of the above results, the condition $\alpha \leq a(x_1) < \beta$ can be rewritten as $\alpha \leq \tilde{c} < \beta$. Now, substituting $a = \tilde{c}$ into (3.3) and solving the resulting equation [noting that $\tilde{r}_-(\tilde{c}) = -\sigma^2\tilde{c}/\mu$], we get

$$V(x) = \frac{M}{\gamma} + \tilde{K}\exp\left(-\frac{\sigma^2\tilde{c}}{\mu}(x - x_1)\right) \qquad \forall x \geq x_1,$$

where $\tilde{K}$ is a constant which takes on either the value of $-\frac{M}{\gamma}$ or $-\mu/(\sigma^2\tilde{c})$ depending, respectively, on whether $x_1 = 0$ or $x_1 > 0$.

Finally, suppose that $a(x_1) < \alpha$. Then it follows from the above results that $\tilde{c} < \alpha$. Therefore $a(x) < \alpha$ for all $x$ in a right neighborhood of $x_1$. Substituting $a = \alpha$ into (3.3) and solving the resulting linear differential equation, we get

$$(3.16) \quad V(x) = \frac{M}{\gamma} + K_1(\alpha)\exp(\tilde{r}_+(\alpha)(x - x_1)) + K_2(\alpha)\exp(\tilde{r}_-(\alpha)(x - x_1)),$$

where $K_1(\alpha)$ and $K_2(\alpha)$ are free constants. If $K_1(\alpha) > 0$, then the right-hand side of (3.16) is unbounded on $[x_1, \infty)$, which contradicts (2.6). If $K_1(\alpha) < 0$, then the right-hand side of (3.16) becomes negative for $x$ large enough, which again is a contradiction. Hence $K_1(\alpha) = 0$. On the other hand, we have $K_2(\alpha) < 0$ in view of $V''(0) < 0$. Therefore,

$$V(x) = \frac{M}{\gamma} + K_\alpha\exp(\tilde{r}_-(\alpha)(x - x_1)),$$

where $K_\alpha$ is a constant that takes on the value either $-\frac{M}{\gamma}$ or $1/(\tilde{r}_-(\alpha))$ depending on whether $x_1$ is zero or not. Combining the above results, we can formulate the following theorem.

THEOREM 3.4. *Let $\tilde{r}_-(\alpha)$, $\tilde{r}_-(\beta)$ and $\tilde{c}$ be the constants given by (3.10) and (3.15), respectively. Let $x_1$ be defined by (3.1). Then for $x_1 = 0$ (resp., for $x_1 > 0$) the following assertions hold.*

(i) *If $\tilde{c} \geq \beta$, then*

$$(3.17) \qquad V(x) = \frac{M}{\gamma} + K_\beta\exp(\tilde{r}_-(\beta)(x - x_1)), \qquad x \geq x_1,$$

*is a concave, twice differentiable solution of the HJB equation (3.3) on $[x_1, \infty)$, where $K_\beta$ is equal to $-\frac{M}{\gamma}$ [resp., to $1/(\tilde{r}_-(\beta))$].*

(ii) *If $\alpha \le \tilde{c} < \beta$, then,*

$$(3.18) \qquad V(x) = \frac{M}{\gamma} + \tilde{K} \exp\left(\frac{-\mu}{\sigma^2 \tilde{c}}(x - x_1)\right), \qquad x \ge x_1,$$

*is a concave, twice differentiable solution of the HJB equation* (3.3) *on* $[x_1, \infty)$, *where $\tilde{K}$ is equal to $-\frac{M}{\gamma}$ [resp., to $-\mu/(\sigma^2 \tilde{c})$].*

(iii) *If $\tilde{c} < \alpha$, then,*

$$(3.19) \qquad V(x) = \frac{M}{\gamma} + K_\alpha \exp(\tilde{r}_-(\alpha)(x - x_1)), \qquad x \ge x_1,$$

*is a concave, twice differentiable solution of the HJB equation* (3.3) *on* $[x_1, \infty)$, *where $K_\alpha$ is a constant equal to $-\frac{M}{\gamma}$ [resp., to $1/(\tilde{r}_-(\alpha))$].*

COROLLARY 3.5. *If $x_1 = 0$, then the solution to* (2.7) *subject to* (2.6) *is given by*

$$V(x) = \begin{cases} \dfrac{M}{\gamma}(1 - \exp(\tilde{r}_-(\beta)x)), & \text{if } \tilde{c} \ge \beta, \\[2mm] \dfrac{M}{\gamma}\left(1 - \exp\left(\dfrac{-\mu}{\sigma^2 \tilde{c}}x\right)\right), & \text{if } \alpha \le \tilde{c} < \beta, \ \forall x \ge 0, \\[2mm] \dfrac{M}{\gamma}(1 - \exp(\tilde{r}_-(\alpha)x)), & \text{if } \tilde{c} < \alpha. \end{cases}$$

Corollary 3.5 shows that the qualitative nature of the solution depends on the relationship between $\frac{2\delta}{\mu}$, $\alpha$ and $\beta$. Accordingly, we consider three cases. However, in contrast to the situation with unbounded dividend rates, each case here will consist of several subcases, each subcase being associated with a different range for the value of $M$.

REMARK 3.6. *If neither $-\frac{M}{\gamma}\tilde{r}_-(\beta) \le 1$ when $\tilde{c} \le \beta$ nor $\frac{M}{\gamma}(\mu/(\sigma^2 \tilde{c})) \le 1$ when $\alpha \le \tilde{c} < \beta$ nor $-\frac{M}{\gamma}\tilde{r}_-(\alpha) \le 1$ when $\tilde{c} < \alpha$ is satisfied, then the solution to* (2.7) *satisfies*

$$V'(0) > 1.$$

The main purpose of the remaining part is to derive the solution to (3.2) and then to combine the latter with Theorem 3.4. The solution to (3.2) is based mainly on the value of $a(0)$. Thus, first of all, we present an analysis of $a(0)$.

PROPOSITION 3.7. *Suppose the assumptions of Remark* 3.6 *hold. Then:*

(i) *$\frac{2\delta}{\mu} < \alpha$ if and only if $a(0) < \alpha$. In this case $a(0) = (\mu\alpha^2)/(2(\mu\alpha - \delta))$.*

(ii) *$\alpha \le \frac{2\delta}{\mu} < \beta$ if and only if $\alpha \le a(0) < \beta$. In this case $a(0) = 2\frac{\delta}{\mu}$.*

(iii) $\beta \le \frac{2\delta}{\mu}$ if and only if $a(0) \ge \beta$. In this case $a(0) = (\mu\beta^2)/(2(\mu\beta - \delta))$.

PROOF. In view of the assumption of Remark 3.6, assume that $x_1$ is positive. Let $\tilde{a} \in [\alpha, \beta]$ such that

$$
\begin{aligned}
(3.20) \qquad 0 &= \max_{\alpha \le a \le \beta} (\tfrac{1}{2}\sigma^2 a^2 V''(0) + (a\mu - \delta)V'(0)) \\
&= \tfrac{1}{2}\sigma^2 \tilde{a}^2 V''(0) + (\tilde{a}\mu - \delta)V'(0).
\end{aligned}
$$

Comparing (3.20) to (3.5) we obtain

$$
(3.21) \qquad \tilde{a}^2 - 2a(0)\tilde{a} + \frac{2\delta}{\mu}a(0) = 0.
$$

From (3.21), it follows that $a(0) \ge \frac{2\delta}{\mu}$. Moreover, by definition, $a(0) \in [\alpha, \ \beta]$ is equivalent to $\tilde{a} = a(0)$, which is further equivalent to $a(0) = \frac{2\delta}{\mu} \in [\alpha, \ \beta]$. Thus we conclude:

(i) If $a(0) < \alpha$, then $\frac{2\delta}{\mu} \le a(0) < \alpha$. Conversely, suppose $\frac{2\delta}{\mu} < \alpha$. If $a(0) \in [\alpha, \beta]$, then by the above results $a(0) = \frac{2\delta}{\mu} < \alpha$, which is a contradiction. Thus either $a(0) < \alpha$ or $a(0) > \beta$. Suppose $a(0) > \beta$. Then $\tilde{a} = \beta$ and by (3.21), $a(0) = (\mu\beta^2)/(2(\mu\beta - \delta)) < \beta$ (due to $\frac{2\delta}{\mu} < \alpha < \beta$). This is again a contradiction. Hence we have $a(0) < \alpha$. Then $\tilde{a} = \alpha$, and in view of (3.21), we get $a(0) = (\mu\alpha^2)/(2(\alpha\mu - \delta))$.

(ii) Suppose $\alpha \le \frac{2\delta}{\mu} < \beta$. Then due to (i) we have $a(0) \ge \alpha$. Now we proceed to prove that $a(0) \le \frac{2\delta}{\mu} < \beta$. Suppose $a(0) > \frac{2\delta}{\mu}$. Then $a(0) > \beta \equiv \tilde{a}$. On the other hand, in view of (3.21), we have $a(0) = (\mu\beta^2)/(2(\beta\mu - \delta))$; thus $(\mu\beta^2)/(2(\beta\mu - \delta)) \ge \beta$, which is equivalent to $2\frac{\delta}{\mu} \ge \beta$. This, however, is a contradiction and therefore $a(0) = \frac{2\delta}{\mu} \in [\alpha, \beta)$. Conversely, if $a(0) \in [\alpha, \beta)$, then $a(0) = \frac{2\delta}{\mu} \in [\alpha, \beta)$.

(iii) Suppose $\beta \le \frac{2\delta}{\mu}$. Then $a(0) \ge \frac{2\delta}{\mu} \ge \beta$, leading to $\tilde{a} = \beta$ and $a(0) = (\mu\beta^2)/(2(\beta\mu - \delta)) \ge \beta$. Conversely, if $a(0) \ge \beta$, then $\tilde{a} = \beta$ and $a(0) = (\mu\beta^2)/(2(\mu\beta - \delta)) \ge \beta$, which is equivalent to $\frac{2\delta}{\mu} \ge \beta$. □

3.1. *Case of $\frac{2\delta}{\mu} < \alpha$.* To resolve equation (3.2), we begin our analysis with an observation that in this case, in view of Proposition 3.7(i), $a(x) < \alpha$ for all $x$ in the right neighborhood of 0. We also suppose that $a(x_1) > \beta$. This assumption is not a restriction, but gives us the solution of (3.2) that corresponds to the maximal interval $[0, \ x_1)$. Substituting $a \equiv \alpha$ in (3.2) and solving the resulting second-order linear ODE, we obtain

$$
(3.22) \qquad V(x) = k_1(\alpha, \beta)(\exp(r_+(\alpha)x) - \exp(r_-(\alpha)x)),
$$

where $k_1(\alpha, \beta)$ is a free constant to be determined and

$$(3.23) \quad \begin{aligned} r_+(z) &= \frac{-(z\mu - \delta) + [(z\mu - \delta)^2 + 2\sigma^2 z^2\gamma]^{1/2}}{\sigma^2 z^2}, \\ r_-(z) &= \frac{-(z\mu - \delta) - [(z\mu - \delta)^2 + 2\sigma^2 z^2\gamma]^{1/2}}{\sigma^2 z^2}, \qquad z > 0. \end{aligned}$$

Due to (3.5) and (3.22),

$$\begin{aligned} a'(x) &= \frac{-\mu}{\sigma^2} \frac{(V''(x))^2 - V'(x)V^{(3)}(x)}{(V''(x))^2} \\ &= \frac{-\mu r_+(\alpha)r_-(\alpha)\exp((r_+(\alpha) + r_-(\alpha))x)(r_+(\alpha) - r_-(\alpha))^2}{\sigma^2(V''(x))^2} > 0 \end{aligned}$$

for each $x$ in the right neighborhood of 0. Therefore $a(x)$ increases and reaches $\alpha$ at the point $x_\alpha$ given by

$$(3.24) \quad x_\alpha = \frac{1}{r_+(\alpha) - r_-(\alpha)} \log\left(\frac{r_-(\alpha)(\mu + \alpha\sigma^2 r_-(\alpha))}{r_+(\alpha)(\mu + \alpha\sigma^2 r_+(\alpha))}\right) > 0.$$

By virtue of Proposition 3.3(i), $\alpha \leq a(x) < \beta$ in the right neighborhood of $x_\alpha$. In this case we substitute $a \equiv a(x)$ and

$$(3.25) \quad V''(x) = \frac{-\mu V'(x)}{\sigma^2 a(x)}$$

into (3.2), differentiating the resulting equation and substituting

$$V''(x) = \frac{-\mu V'(x)}{\sigma^2 a(x)}$$

once more, we arrive at

$$\frac{\mu a'(x)}{2} + \frac{\mu\delta}{\sigma^2 a(x)} = \frac{\mu^2 + 2\gamma\sigma^2}{2\sigma^2}.$$

As a result,

$$(3.26) \quad a'(x) = \frac{\mu^2 + 2\gamma\sigma^2}{\mu\sigma^2}\left(1 - \frac{c}{a(x)}\right)$$

with

$$(3.27) \quad c \equiv \frac{2\delta\mu}{\mu^2 + 2\gamma\sigma^2}.$$

Integrating (3.26), we get $G(a(x)) \equiv ((\mu^2 + 2\gamma\sigma^2)/(\mu\sigma^2))(x - x_\alpha) + G(\alpha)$, where

$$(3.28) \quad G(u) = u + c\log(u - c).$$

Therefore,

$$(3.29) \qquad a(x) = G^{-1}\left(\frac{\mu^2 + 2\gamma\sigma^2}{\mu\sigma^2}(x - x_\alpha) + G(\alpha)\right).$$

Thus $a(x)$ is increasing and $a(x_\beta) = \beta$ for

$$(3.30) \qquad \begin{aligned} x_\beta &\equiv \frac{\mu\sigma^2}{\mu^2 + 2\gamma\sigma^2}[G(\beta) - G(\alpha)] + x_\alpha \\ &= \frac{\mu\sigma^2}{\mu^2 + 2\gamma\sigma^2}(\beta - \alpha) + \frac{\mu\sigma^2 c}{\mu^2 + 2\gamma\sigma^2}\log\left(\frac{\beta - c}{\alpha - c}\right). \end{aligned}$$

Solving (3.25) we obtain

$$(3.31) \quad V(x) = V(x_\alpha) + V'(x_\alpha)\int_{x_\alpha}^x \exp\left(-\frac{\mu}{\sigma^2}\int_{x_\alpha}^y \frac{du}{a(u)}\right) dy, \qquad x_\alpha \le x < x_\beta,$$

where $V(x_\alpha)$ and $V'(x_\alpha)$ are free constants. Choosing $V(x_\alpha)$ and $V'(x_\alpha)$ as the value and the derivative, respectively, of the right-hand side of (3.22) at $x_\alpha$, we can ensure that the function $V$ given by (3.22) and (3.31) is continuous with its first and second derivatives at the point $x_\alpha$ no matter what the choice of $k(\alpha, \beta)$ is. (Note that due to the HJB equation, continuity of $V$ and its first derivative at $x_\alpha$ automatically implies continuity of the second derivative as well.) Next we simplify (3.31). First, changing variables $a(u) = \theta$ we get

$$\int_{x_\alpha}^x \exp\left(-\frac{\mu}{\sigma^2}\int_{x_\alpha}^y \frac{du}{a(u)}\, dy\right)$$
$$= \frac{\mu\sigma^2}{\mu^2 + 2\gamma\sigma^2}\int_\alpha^{a(x)}\left(1 + \frac{c}{\theta - c}\right)\left(\frac{\theta - c}{\alpha - c}\right)^{-\Gamma} d\theta, \qquad x_\alpha \le x < x_\beta.$$

On the other hand, relationships (3.24) and (3.22) imply

$$V(x_\alpha) = \frac{\alpha\mu - 2\delta}{2\gamma}V'(x_\alpha).$$

Simple algebraic transformations yield

$$(3.32) \qquad \left(\frac{\mu\sigma^2}{\mu^2 + 2\gamma\sigma^2}\right)\left(\frac{c}{\Gamma} - \frac{z - c}{1 - \Gamma}\right) = \frac{z\mu - 2\delta}{2\gamma} \qquad \forall z > 0,$$

where $c$ is given by (3.27) and

$$(3.33) \qquad \Gamma = \frac{\mu^2}{\mu^2 + 2\gamma\sigma^2}.$$

Therefore,

$$(3.34) \qquad V(x) = V'(x_\alpha)\frac{\mu a(x) - 2\delta}{2\gamma}\left(\frac{a(x) - c}{\alpha - c}\right)^{-\Gamma}, \qquad x_\alpha \le x < x_\beta.$$

The same arguments as in Proposition 3.3(ii) show that $a(x) \geq \beta$ for each $x \geq x_\beta$. Thus, substituting $a \equiv \beta$ into (3.2) and solving the resulting ODE, we get

$$(3.35) \quad V(x) = k_1(\beta)\exp(r_+(\beta)(x-x_1)) + k_2(\beta)\exp(r_-(\beta)(x-x_1)),$$
$$x_\beta \leq x < x_1,$$

where $k_1(\beta)$ and $k_2(\beta)$ are two free constants to be determined. The continuity of (3.31) at $x_\beta$, together with simple but tedious algebraic transformation [similar to those used above to simplify (3.31) to (3.34)] lead to

$$(3.36) \quad V'(x_\alpha) = V'(x_\beta)\left(\frac{\beta - c}{\alpha - c}\right)^\Gamma.$$

Let $\tilde{c}$ be given by (3.15) and

$$(3.37) \quad M_z = \left(z - \frac{2\delta\mu}{\mu^2 + 2\sigma^2\gamma}\right)\frac{\mu^2 + 2\sigma^2\gamma}{2\mu}, \qquad z > 0.$$

Then $\tilde{c} \geq \beta$ (resp., $\tilde{c} = \beta$ ) is equivalent to $M \geq M_\beta$ (resp., $M = M_\beta$). This is the first subcase we consider.

3.1.1. *Case of $M > M_\beta$.* Our assumptions imply that in this case, $a(x_1) > \beta$, which is equivalent to $x_1 > x_\beta$. Combining (3.22), (3.34) and (3.35) and Theorem 3.4(i) we can write a general form of the solution to (2.7) and (2.6),

$$(3.38) \ V(x) = \begin{cases} K_1(\alpha, \beta)(\exp(r_+(\alpha)x) - \exp(r_-(\alpha)x)), & 0 \leq x < x_\alpha, \\ V'(x_\alpha)\dfrac{\mu a(x) - 2\delta}{2\gamma}\left(\dfrac{a(x) - c}{\alpha - c}\right)^{-\Gamma}, & x_\alpha \leq x < x_\beta, \\ K_1(\beta)\exp(r_+(\beta)(x-x_1)) \\ \quad + K_2(\beta)\exp(r_-(\beta)(x-x_1)), & x_\beta \leq x < x_1, \\ \dfrac{M}{\gamma} + \dfrac{1}{\tilde{r}_-(\beta)}\exp(\tilde{r}_-(\beta)(x-x_1)), & x \geq x_1, \end{cases}$$

where $r_+(\alpha)$, $r_-(\alpha)$, $r_+(\beta)$ and $r_-(\beta)$, $x_\alpha$ and $x_\beta$ are given by (3.23), (3.24) and (3.30), respectively, and $K_1(\beta)$, $K_2(\beta)$, $K_1(\alpha, \beta)$ and $x_1$ are unknown constants to be determined. Continuity of the first and the second derivatives at $x_1$ results in

$$V'(x_1) = 1, \qquad V''(x_1) = \tilde{r}_-(\beta).$$

This gives us two equations,

$$1 = K_1(\beta)r_+(\beta) + K_2(\beta)r_-(\beta),$$
$$\tilde{r}_-(\beta) = K_1(\beta)r_+^2(\beta) + K_2(\beta)r_-^2(\beta),$$

whose solutions are

$$(3.39) \qquad K_1(\beta) = \frac{\tilde{r}_-(\beta) - r_-(\beta)}{r_+(\beta)(r_+(\beta) - r_-(\beta))},$$

$$K_2(\beta) = \frac{r_+(\beta) - \tilde{r}_-(\beta)}{r_-(\beta)(r_+(\beta) - r_-(\beta))}.$$

Put $\Delta = x_\beta - x_1$. As before, using the principle of smooth fit at $x_\beta$, we get

$$(3.40) \qquad x_\beta - x_1 = \Delta = \frac{1}{(r_+(\beta) - r_-(\beta))}$$

$$\times \log\left(-\frac{(r_+(\beta) - \tilde{r}_-(\beta))(\mu + \beta\sigma^2 r_-(\beta))}{(\tilde{r}_-(\beta) - r_-(\beta))(\mu + \beta\sigma^2 r_+(\beta))}\right).$$

The expression on the right-hand side of (3.40) is negative due to $\tilde{c} > \beta$. In view of (3.40) and (3.39) we can derive a simplified expression for $V'(x_\beta)$:

$$V'(x_\beta) = \beta\sigma^2 \exp(r_-(\beta)\Delta)\frac{r_+(\beta) - \tilde{r}_-(\beta)}{\mu + \beta\sigma^2 r_+(\beta)}.$$

The continuity of $V'$ at $x_\alpha$ yields $V'(x_\alpha) = K_1(\alpha, \beta)(r_+(\alpha)\exp(r_+(\alpha)x_\alpha) - r_-(\alpha)\exp(r_-(\alpha)x_\alpha))$. Combining this equality with (3.36), we get

$$(3.41) \qquad K_1(\alpha, \beta) = \frac{V'(x_\beta)((\beta - c)/(\alpha - c))^\Gamma}{r_+(\alpha)\exp(r_+(\alpha)x_\alpha) - r_-(\alpha)\exp(r_-(\alpha)x_\alpha)}.$$

THEOREM 3.8. *Let $a(x)$ be a function given by (3.29) and let $r_+(\alpha)$, $r_-(\alpha)$, $r_+(\beta)$, $r_-(\beta)$, $x_\alpha$, $x_\beta$, $c$, $\Gamma$, $\tilde{r}_-(\beta)$, $K_1(\beta)$, $K_2(\beta)$, $K_1(\alpha, \beta)$ and $x_1$ be given by (3.23), (3.24), (3.30), (3.27), (3.33), (3.10), (3.39), (3.41) and (3.40), respectively. If $\frac{2\delta}{\mu} < \alpha$ and $M > M_\beta$, then $V$ given by (3.38) is a concave, twice differentiable solution of the HJB equation (2.7), subject to (2.6).*

PROOF. From the way we constructed $V$, it is a twice continuously differentiable solution to the HJB equation (3.2). What remains to show is the concavity. From (3.38), we deduce that

$$V'''(x) = k_1(\alpha, \beta)(r_+^3(\alpha)\exp(r_+(\alpha)x) - r_-^3(\alpha)\exp(r_-(\alpha)x)) > 0$$

$$\forall 0 \leq x < x_\alpha,$$

due to $r_-(\alpha) < 0 < k_1(\alpha, \beta)$. Hence on this interval $V''$ is increasing and

$$V''(x) < V''(x_\alpha) = k_1(\alpha, \beta)(r_+^2(\alpha)\exp(r_+(\alpha)x_\alpha) - r_-^2(\alpha)\exp(r_-(\alpha)x_\alpha)) < 0,$$

due to $(r_-(\alpha))/(r_+(\alpha)) = \exp((r_+(\alpha) - r_-(\alpha))x_\alpha)$ and $|r_-(\alpha)| > r_+(\alpha)$.

For $x_\alpha \le x < x_\beta$, $V''(x) = (-\mu V'(x))/(\sigma^2 a(x)) < 0$. For $x_\beta \le x < x_1$,

$$V'''(x) = k_1(\beta)r_+^3(\beta)\exp(r_+(\beta)(x - x_1)) + k_2(\beta)r_-^3(\beta)\exp(r_-(\beta)(x - x_1)) > 0,$$

since $k_2(\beta)$ and $r_-(\beta)$ are of the same sign. Thus $V''(x) < V''(x_1) < 0$ $\forall x_\beta \le x < x_1$. Finally, $V''(x) < 0$ $\forall x \ge x_1$. This establishes the concavity of $V$. Since $V'(x) > 1$ for $x < x_1$ and $V'(x) \le 1$ for $x \ge x_1$, it is clear that $V$ satisfies (2.7). $\qquad \square$

3.1.2. *Case of $M_\alpha < M \le M_\beta$.* Expression (3.37) shows that $\beta \ge \tilde{c} > \alpha$ if and only if $M_\beta \ge M > M_\alpha$. From (3.29) we see that the condition $\beta \ge \tilde{c} > \alpha$ is equivalent to $\beta \ge a(x_1) > \alpha$. In view of Theorem 3.4, $a(x) \le a(x_1) \le \beta$ for all $x \ge 0$. This also implies $V'(x_1) = 1$. As a result $\tilde{K} = -\sigma^2 \tilde{c}/\mu$. Taking into account (3.22), (3.34) and (3.18), we can write the expression for $V$ as

$$V(x) = \begin{cases} K_1(\alpha,\beta)(\exp(r_+(\alpha)x) - \exp(r_-(\alpha)x)), & 0 \le x < x_\alpha, \\ V'(x_\alpha)\dfrac{\mu a(x) - 2\delta}{2\gamma}\left(\dfrac{a(x) - c}{2\delta/\mu - c}\right)^{-\Gamma}, & x_\alpha \le x < x_1, \\ \dfrac{M}{\gamma} - \dfrac{\sigma^2 \tilde{c}}{\mu}\exp\left(-\dfrac{\mu}{\sigma^2 \tilde{c}}(x - x_1)\right), & x \ge x_1. \end{cases}$$

For $a(x)$ given by (3.29), the root of the equation $a(x_1) = \tilde{c}$ can be written as

$$(3.42) \qquad \begin{aligned} x_1 &= \frac{\mu\sigma^2}{\mu^2 + 2\sigma^2\gamma}\int_{2\delta/\mu}^{\tilde{c}} \frac{u\,du}{u - c} \\ &= \frac{\mu\sigma^2(\tilde{c} - 2\delta/\mu)}{\mu^2 + 2\sigma^2\gamma} + \frac{\mu\sigma^2}{\mu^2 + 2\sigma^2\gamma}\log\left(\frac{\tilde{c} - c}{2\delta/\mu - c}\right). \end{aligned}$$

The continuity of $V'$ at $x_1$ leads to $V'(x_\alpha) = (\frac{\tilde{c} - c}{2\delta/\mu - c})^\Gamma$. Consequently,

$$(3.43) \qquad K_1(\alpha,\beta) = \frac{((\tilde{c} - c)/(2\delta/\mu - c))^\Gamma}{r_+(\alpha)\exp(r_+(\alpha)x_\alpha) - r_-(\alpha)\exp(r_-(\alpha)x_\alpha)}.$$

THEOREM 3.9. *Let $a(x)$ be a function given by (3.29) and let $r_+(\alpha)$, $r_-(\alpha)$, $K_1(\alpha,\beta)$, $x_\alpha$, $x_1$, $c$, $\Gamma$ and $\tilde{c}$ be given by (3.23), (3.43), (3.24), (3.42), (3.27), (3.33) and (3.15), respectively. If $\frac{2\delta}{\mu} < \alpha$ and $M_\alpha < M \le M_\beta$, then*

$$(3.44)\ V(x) = \begin{cases} K_1(\alpha,\beta)(\exp(r_+(\alpha)x) - \exp(r_-(\alpha)x)), & 0 \le x < x_\alpha, \\ \dfrac{\mu a(x) - 2\delta}{2\gamma}\left(\dfrac{a(x) - c}{\tilde{c} - c}\right)^{-\Gamma}, & x_\alpha \le x < x_1, \\ \dfrac{M}{\gamma} - \dfrac{\sigma^2 \tilde{c}}{\mu}\exp\left(-\dfrac{\mu}{\sigma^2 \tilde{c}}(x - x_1)\right), & x \ge x_1 \end{cases}$$

*is a concave, twice differentiable solution of the HJB equation (2.7) subject to (2.6).*

PROOF. The proof of this theorem follows the same lines as that of Theorem 3.8. □

Now suppose that $M \leq M_\alpha$. Then $a(x) \leq a(x_1) \leq \alpha$ for each $x \geq 0$ [since $a(x)$ is increasing on $[0, x_1)$ and is constant for $x \geq x_1$; see Theorem 3.4]. If $V'(0) > 1$, then $x_1 > 0$ and $V'(x_1) = 1$. As a result, $\tilde{K}_\alpha = 1/(\tilde{r}_-(\alpha))$. In view of (3.22) and (3.19), the function $V$ is given by

$$(3.45) \quad V(x) = \begin{cases} k_1(\alpha, \beta)(\exp(r_+(\alpha)x) - \exp(r_-(\alpha)x)), & 0 \leq x < x_1, \\ \dfrac{M}{\gamma} + \dfrac{1}{\tilde{r}_-(\alpha)} \exp(\tilde{r}_-(\alpha)(x - x_1)), & x \geq x_1. \end{cases}$$

The smoothness of $V$ requires

$$V'(x_1-) = 1, \qquad V''(x_1-) = \tilde{r}_-(\alpha),$$

which translates into

$$(3.46) \quad \begin{aligned} k_1(\alpha, \beta)(r_+(\alpha)\exp(r_+(\alpha)x_1) - r_-(\alpha)\exp(r_-(\alpha)x_1)) &= 1, \\ k_1(\alpha, \beta)(r_+^2(\alpha)\exp(r_+(\alpha)x_1) - r_-^2(\alpha)\exp(r_-(\alpha)x_1)) &= \tilde{r}_-(\alpha). \end{aligned}$$

Excluding $k_1(\alpha, \beta)$, we get an equation for $x_1$:

$$(3.47) \quad \exp((r_+(\alpha) - r_-(\alpha))x_1) = \frac{r_-(\alpha)(r_-(\alpha) - \tilde{r}_-(\alpha))}{r_+(\alpha)(r_+(\alpha) - \tilde{r}_-(\alpha))}.$$

This equation has a positive solution if and only if

$$(3.48) \quad M > M_0(\alpha) \equiv \frac{\alpha^2 \sigma^2 \gamma}{2(\alpha\mu - \delta)}.$$

This proves the following statement.

PROPOSITION 3.10. *If $\frac{2\delta}{\mu} < \alpha$, then*

$$V'(0) > 1 \qquad iff \ M > \frac{\alpha^2 \sigma^2 \gamma}{2(\alpha\mu - \delta)}.$$

Let $M_z$ be given by (3.37) and

$$M_0(z) \equiv \frac{z^2 \sigma^2 \gamma}{2(z\mu - \delta)}, \qquad z > \frac{\delta}{\mu}.$$

A simple analysis shows that $f(z) \equiv M_0(z) - M_z$ is a decreasing function of $z$ and $f(\frac{2\delta}{\mu}) = 0$. Similarly, we claim that $M_0(z)$ is decreasing for $z \leq \frac{2\delta}{\mu}$ and increasing for $z \geq \frac{2\delta}{\mu}$. Thus, we derive the inequalities

$$(3.49) \qquad M_0\left(\frac{2\delta}{\mu}\right) < M_0(\alpha) < M_\alpha < M_\beta \qquad \text{if } \frac{2\delta}{\mu} < \alpha,$$

$$(3.50) \qquad M_\alpha \leq M_0(\alpha) \leq M_0\left(\frac{2\delta}{\mu}\right) < M_\beta \qquad \text{if } \alpha \leq \frac{2\delta}{\mu} < \beta,$$

$$(3.51) \qquad M_\alpha < M_\beta \leq M_0\left(\frac{2\delta}{\mu}\right) \leq M_0(\beta) < M_0(\alpha) \qquad \text{if } \beta \leq \frac{2\delta}{\mu}.$$

Since the qualitative behavior of the solution to (3.2) [resp., to (3.3)] depends on the value of $a(0)$ [resp., of $a(x_1)$], in accordance with (3.49) we distinguish and study the remaining subcases in the following sections.

3.1.3. *Case of $M_0(\alpha) < M \leq M_\alpha$.* This is the case when (3.47) has a positive solution $x_1$ given by

$$(3.52) \qquad x_1 = \frac{1}{r_+(\alpha) - r_-(\alpha)} \log\left(\frac{r_-(\alpha)(r_-(\alpha) - \tilde{r}_-(\alpha))}{r_+(\alpha)(r_+(\alpha) - \tilde{r}_-(\alpha))}\right).$$

THEOREM 3.11. *Let $r_+(\alpha)$, $r_-(\alpha)$, $\tilde{r}_-(\alpha)$ and $x_1$ be given by (3.23), (3.10) and (3.52), respectively, and let $k_1(\alpha, \beta)$ be determined by (3.46). If $\frac{2\delta}{\mu} < \alpha$ and $M_0(\alpha) < M \leq M_\alpha$, then $V$ given by (3.45) is a concave, twice continuously differentiable solution of (2.7) subject to (2.6).*

PROOF. The proof of this theorem follows from that of Theorem 3.9. □

3.1.4. *Case of $M \leq M_0(\alpha)$.* By virtue of Proposition 3.10, this assumption results in $V'(0) \leq 1$. As a consequence, $x_1 = 0$. As shown in Theorem 3.4, this leads to $a(x) = a(0)$ for each $x \geq 0$. Since $M_\alpha > M_0(\alpha)$, we can apply Corollary 3.5 to deduce $a(0) < \alpha$.

THEOREM 3.12. *Let $\tilde{r}_-(\alpha)$ be a constant given by (3.10). If $\frac{2\delta}{\mu} < \alpha$ and $M \leq M_0(\alpha)$, then*

$$(3.53) \qquad V(x) = \frac{M}{\gamma}(1 - \exp(\tilde{r}_-(\alpha)x)), \qquad x \geq 0,$$

*is a concave, twice continuously differentiable solution of (2.7) subject to (2.6).*

PROOF. See Corollary 3.5. □

3.2. *Case of $\alpha \leq \frac{2\delta}{\mu} < \beta$.* In this section, we investigate the second main case of $\alpha \leq a(0) < \beta$. As in the preceding section, if we assume $a(x_1) > \beta$, then (3.2) admits the solution

$$(3.54) \quad V(x) = \begin{cases} V'(0)\dfrac{\mu a(x) - 2\delta}{2\gamma}\left(\dfrac{a(x) - c}{2\delta/\mu - c}\right)^{-\Gamma}, & 0 \leq x < x_\beta, \\ K_1(\beta)\exp(r_+(\beta)(x - x_1)) \\ \quad + K_2(\beta)\exp(r_-(\beta)(x - x_1)), & x_\beta \leq x < x_1. \end{cases}$$

Here

$$(3.55) \quad \begin{aligned} x_\beta &= \frac{\mu\sigma^2}{\mu^2 + 2\gamma\sigma^2}\left[G(\beta) - G\left(\frac{2\delta}{\mu}\right)\right] \\ &= \frac{\mu\sigma^2}{\mu^2 + 2\gamma\sigma^2}\left(\beta - \frac{2\delta}{\mu}\right) + \frac{2\delta\mu c}{\mu^2 + 2\gamma\sigma^2}\log\left(\frac{\beta - c}{2\delta/\mu - c}\right) \end{aligned}$$

and the function $a(x)$ is defined by

$$(3.56) \quad a(x) = G^{-1}\left(\frac{\mu^2 + 2\gamma\sigma^2}{\mu\sigma^2}x + G\left(\frac{2\delta}{\mu}\right)\right) \in \left[\frac{2\delta}{\mu}, \infty\right),$$

where $G$ is given by (3.28).

As before, the solution to (2.7) is derived by combining (3.54) and (3.3), by distinguishing subcases as follows.

3.2.1. *Case of $M > M_\beta$.* As in Section 3.1.1, consider the case of $x_1 > x_\beta$. This case is characterized by $M > M_\beta$, which is also equivalent to $\tilde{c} > \beta$. As a result, we get $V'(x_1) = 1$. This leads to $\tilde{K}_\beta = 1/(\tilde{r}_-(\beta))$ [see Theorem 3.4(i)]. Then using (3.54) and Theorem 3.4(i), $V$ can be represented in the form

$$(3.57) \quad V(x) = \begin{cases} V'(0)\dfrac{\mu a(x) - 2\delta}{2\gamma}\left(\dfrac{a(x) - c}{2\delta/\mu - c}\right)^{-\Gamma}, & 0 \leq x < x_\beta, \\ K_1(\beta)\exp(r_+(\beta)(x - x_1)) \\ \quad + K_2(\beta)\exp(r_-(\beta)(x - x_1)), & x_\beta \leq x < x_1, \\ \dfrac{M}{\gamma} + \dfrac{1}{\tilde{r}_-(\beta)}\exp(\tilde{r}_-(\beta)(x - x_1)), & x \geq x_1. \end{cases}$$

where $a(x)$, $K_1(\beta)$, $K_2(\beta)$ and $x_1$ are given by (3.56), (3.39) and (3.40), respectively. Let $\Delta = x_\beta - x_1$ be given by (3.30). Continuity of $V$ at $x_\beta$ yields

$$(3.58) \quad \begin{aligned} V'(0) &= \left(\frac{2\gamma}{\mu\beta - 2\delta}\left(\frac{\beta - c}{2\delta/\mu - c}\right)^\Gamma\right) \\ &\quad \times (K_1(\beta)\exp(r_+(\beta)\Delta) + K_2(\beta)\exp(r_-(\beta)\Delta)). \end{aligned}$$

THEOREM 3.13.  *Let $V'(0)$, $a(x)$, $c$, $\Gamma$, $K_1(\beta)$, $K_2(\beta)$, $x_1$ and $\tilde{c}$ be given by* (3.58), (3.56), (3.27), (3.33), (3.39), (3.40) *and* (3.15), *respectively. If $\alpha \leq \frac{2\delta}{\mu} < \beta$ and $M \geq M_\beta$, then $V(x)$ given by* (3.57) *is a concave, twice continuously differentiable solution of* (2.7) *subject to* (2.6).

PROOF.  The proof results from combining (3.54) and Theorem 3.4(ii). □

To classify the remaining cases, suppose $M \leq M_\beta$. Let $V'(0) > 1$. In this case, $x_1$ defined by (3.1) is positive. Therefore, using (3.54) and Theorem 3.4(ii), we can represent $V$ as

$$(3.59) \quad V(x) = \begin{cases} V'(0)\dfrac{\mu a(x) - 2\delta}{2\gamma}\left(\dfrac{a(x) - c}{2\delta/\mu - c}\right)^{-\Gamma}, & 0 \leq x < x_1, \\ V(x_1) + \dfrac{\sigma^2 \tilde{c}}{\mu} - \dfrac{\sigma^2 \tilde{c}}{\mu}\exp\left(-\dfrac{\mu}{\sigma^2 \tilde{c}}(x - x_1)\right), & x \geq x_1, \end{cases}$$

where $a(x)$ and $\tilde{c}$ are given by (3.56) and (3.15), respectively. As a consequence, we get

$$(3.60) \quad a(x_1) = \tilde{c}.$$

Continuity of $V'(x)$ at $x = x_1$ [see (3.26) for the expression of the derivatives of $a(x)$] along with (3.60) results in

$$(3.61) \quad V'(0) = \left(\frac{\tilde{c} - c}{2\delta/\mu - c}\right)^{\Gamma}.$$

Substituting (3.61), (3.60) and (3.37) into (3.59), we obtain

$$(3.62) \quad V(x_1) = \frac{M}{\gamma} - \frac{\sigma^2 \tilde{c}}{\mu}.$$

The unknown constant $x_1$ is the root of (3.60). Recalling (3.56), we see that (3.60) admits a positive solution if and only if

$$(3.63) \quad M > M_0\left(\frac{2\delta}{\mu}\right) = \frac{2\delta^2 \sigma^2 \gamma}{\mu^2}.$$

PROPOSITION 3.14.  *Suppose $\alpha \leq \frac{2\delta}{\mu} < \beta$. Then*

$$V'(0) > 1 \qquad iff\ M > M_0\left(\frac{2\delta}{\mu}\right).$$

In view of this proposition, we distinguish the remaining subcases as follows.

3.2.2. *Case of $M_0(\frac{2\delta}{\mu}) < M \le M_\beta$.* Substituting (3.56) into (3.60), we obtain

$$
\begin{aligned}
x_1 &= \frac{\mu\sigma^2}{\mu^2 + 2\sigma^2\gamma} \int_{2\delta/\mu}^{\tilde{c}} \frac{u\,du}{u-c} \\
&= \frac{\mu\sigma^2(\tilde{c} - 2\delta/\mu)}{\mu^2 + 2\sigma^2\gamma} + \frac{\mu\sigma^2}{\mu^2 + 2\sigma^2\gamma} \log\left(\frac{\tilde{c} - c}{2\delta/\mu - c}\right).
\end{aligned}
\tag{3.64}
$$

THEOREM 3.15. *Let $a(x)$, $c$, $\Gamma$, $\tilde{c}$ and $x_1$ be given by* (3.56), (3.27), (3.33), (3.15) *and* (3.64), *respectively. If $\alpha \le \frac{2\delta}{\mu} < \beta$ and $M_0(\frac{2\delta}{\mu}) < M \le M_\beta$, then*

$$
V(x) = \begin{cases} \dfrac{\mu a(x) - 2\delta}{2\gamma}\left(\dfrac{a(x) - c}{\tilde{c} - c}\right)^{-\Gamma}, & 0 \le x < x_1, \\[2ex] \dfrac{M}{\gamma} - \dfrac{\sigma^2\tilde{c}}{\mu}\exp\left(-\dfrac{\mu}{\sigma^2\tilde{c}}(x - x_1)\right), & x \ge x_1, \end{cases}
\tag{3.65}
$$

*is a concave, twice continuously differentiable solution of* (2.7) *subject to* (2.6).

PROOF. The proof of this theorem follows from that of (3.54). □

3.2.3. *Case of $M \le M_0(\frac{2\delta}{\mu})$.* Note that in this case, $V'(0) \le 1$ due to Proposition 3.14. From (3.50), it follows that $a(0) = \tilde{c} < \beta$.

THEOREM 3.16. *Suppose $\alpha \le \frac{2\delta}{\mu} < b$ and $M \le M_0(\frac{2\delta}{\mu})$. Let $\tilde{c}$ and $\tilde{r}_-(\alpha)$ be given by* (3.15) *and* (3.10), *respectively.*

(i) *If $\alpha < \frac{2\delta}{\mu}$ and $\alpha \le \tilde{c}$, then*

$$
V(x) = \frac{M}{\gamma}\left(1 - \exp\left(-\frac{\mu}{\sigma^2\tilde{c}}x\right)\right), \qquad x \ge 0,
\tag{3.66}
$$

*is a concave, twice continuously differentiable solution of* (2.7) *subject to* (2.6). *If $\alpha = \frac{2\delta}{\mu}$ or $\tilde{c} < \alpha$, then*

$$
V(x) = \frac{M}{\gamma}(1 - \exp(\tilde{r}_-(\alpha)x)), \qquad x \ge 0,
\tag{3.67}
$$

*is a concave, twice continuously differentiable solution of* (2.7) *subject to* (2.6).

PROOF. See Corollary 3.5. □

3.3. *Case of $\frac{\delta}{\mu} < \beta \leq \frac{2\delta}{\mu}$.* We now investigate the final main case, $\frac{\delta}{\mu} < \beta \leq \frac{2\delta}{\mu}$. Suppose $V'(0) > 1$. Then $x_1$ defined by (3.1) is positive. In view of Proposition 3.7(iii) and Proposition 3.3(i), $a(x) \geq \beta$ for all $x \geq 0$. Therefore, using (3.62) and Theorem 3.4(i), $V$ can be represented as

$$(3.68) \quad V(x) = \begin{cases} K(\exp(r_+(\beta)x) - \exp(r_-(\beta)x)), & 0 \leq x < x_1, \\ \dfrac{M}{\gamma} + \dfrac{1}{\tilde{r}_-(\beta)} \exp(\tilde{r}_-(\beta)(x - x_1)), & x \geq x_1. \end{cases}$$

The principle of smooth fit for $V$ at $x_1$ yields

$$(3.69) \quad \begin{aligned} & V'(x_1-) = K(r_+(\beta)\exp(r_+(\beta)x_1) - r_-(\beta)\exp(r_-(\beta)x_1)) = 1, \\ & V''(x_1-) = \tilde{r}_-(\beta). \end{aligned}$$

Thus

$$(3.70) \quad \exp((r_+(\beta) - r_-(\beta))x_1) = \frac{r_-(\beta)(r_-(\beta) - \tilde{r}_-(\beta))}{r_+(\beta)(r_+(\beta) - \tilde{r}_-(\beta))},$$

which admits a positive solution $x_1$ iff

$$(3.71) \quad M > M_0(\beta) = \frac{\sigma^2 \beta^2 \gamma}{2(\beta\mu - \delta)}.$$

PROPOSITION 3.17. *If $\frac{\delta}{\mu} < \beta \leq \frac{2\delta}{\mu}$, then*

$$V'(0) > 1 \qquad \textit{iff } M > M_0(\beta).$$

PROOF. The proof of this proposition follows from the calculations in this and the previous sections. □

In view of Proposition 3.17, we need to treat only two subcases, namely, $M > M_0(\beta)$ and $M \leq M_0(\beta)$, to complete our analysis.

3.3.1. *Case of $M > M_0(\beta)$.* In this case, (3.70) has a positive solution $x_1$ given by

$$(3.72) \quad x_1 = \frac{1}{r_+(\beta) - r_-(\beta)} \log\left(\frac{r_-(\beta)(r_-(\beta) - \tilde{r}_-(\beta))}{r_+(\beta)(r_+(\beta) - \tilde{r}_-(\beta))}\right).$$

THEOREM 3.18. *Let $x_1$, $r_+(\beta)$, and $r_-(\beta)$ and $\tilde{r}_-(\beta)$ be given by (3.72), (3.23) and (3.10), respectively, and let $K$ be a constant determined from (3.69). If $\frac{\delta}{\mu} < \beta \leq \frac{2\delta}{\mu}$ and $M > M_0(\beta)$, then $V$ given by (3.68) is a concave, twice continuously differentiable solution of (2.7) subject to (2.6).*

PROOF. By differentiating the expression (3.68), we obtain

$$V^{(3)}(x) = K(r_+^3(\beta)\exp(r_+(\beta)x) - r_-^3(\beta)\exp(r_-(\beta)x)) > 0, \qquad 0 \leq x < x_1.$$

As a result, $V''(x) < V''(x_1) = 0$ and $V'(x) > V'(x_1) = 1$. This proves that $V$ is concave. $\square$

3.3.2. *Case of $M \leq M_0(\beta)$.* In this case, in view of Proposition 3.17, $V'(0) \leq 1$. Therefore, $x_1$ defined by (3.1) equals zero.

THEOREM 3.19. *Suppose that either $\beta \leq \frac{\delta}{\mu}$ or $\frac{\delta}{\mu} < \beta \leq \frac{2\delta}{\mu}$ and $M \leq M_0(\beta)$. Let $\tilde{r}_-(\beta)$, $\tilde{r}_-(\alpha)$ and $\tilde{c}$ be given by (3.10) and (3.15), respectively.*

(i) *If $M \geq M_\beta$, then*

$$(3.73) \qquad V(x) = \frac{M}{\gamma}(1 - \exp(\tilde{r}_-(\beta)x)), \qquad x \geq 0,$$

*is a concave, twice continuously differentiable solution of (2.7) subject to (2.6).*

(ii) *If $M_\alpha \leq M < M_\beta$, then*

$$V(x) = \frac{M}{\gamma}\left(1 - \exp\left(-\frac{\mu}{\sigma^2 \tilde{c}}x\right)\right), \qquad x \geq 0,$$

*is a concave, twice continuously differentiable solution of (2.7) subject to (2.6).*

(iii) *If $M < M_\alpha$, then*

$$V(x) = \frac{M}{\gamma}(1 - \exp(\tilde{r}_-(\alpha)x)), \qquad x \geq 0,$$

*is a concave, twice continuously differentiable solution of (2.7) subject to (2.6).*

PROOF. In the case of $\beta \leq \frac{\delta}{\mu}$, the inequality $V'(0) \leq 1$ holds due to Proposition 3.17. Then by applying Corollary 3.5, the desired result follows. On the other hand, if $\frac{\delta}{\mu} < \beta \leq \frac{2\delta}{\mu}$ and $M \leq M_0(\beta)$, then $V'(0) \leq 1$ (see Proposition 3.1). Thus, in view of (3.51), Corollary 3.5 can be applied again to obtain the results. $\square$

## 4. Optimal policies.

In this section we construct the optimal control policies based on the solutions to the HJB equations obtained in the previous section. The derivation of the results of this section is simpler than the corresponding one in [4], in view of the fact that no Skorohod problem has to be involved in this case.

Suppose $V$ is a concave solution to the HJB equation (2.7). Define

$$
\begin{aligned}
(4.1) \quad a^*(x) \equiv \arg\max_{\alpha \le a \le \beta} (\tfrac{1}{2}\sigma^2 a^2 V''(x) + (a\mu - \delta)V'(x) \\
- \gamma V(x) + M(1 - V'(x))^+)
\end{aligned}
$$

and

$$
M^*(x) \equiv M\mathbb{1}_{\{x \ge x_1\}},
$$

where $x_1$ is defined by (3.1). The function $a^*(x)$ is the optimal feedback risk control function, while the function $M^*(x)$ represents the optimal dividend rate payments when the level of the reserve is $x$.

THEOREM 4.1. *Let $R_t^*$, $t \ge 0$, be a solution to the stochastic differential equation*

$$
\begin{aligned}
(4.2) \quad dR_t^* &= [a^*(R_t^*)\mu - \delta - M^*(R_t^*)]\,dt + a^*(R_t^*)\sigma\,dW_t, \\
R_0^* &= x.
\end{aligned}
$$

*Then for $\pi^* \equiv (a_t^*, c_t^*; t \ge 0) = (a^*(R_t^*), M^*(R_t^*); t \ge 0)$, we have*

$$
(4.3) \quad J_x(\pi^*) = V(x) \qquad \forall\, x \ge 0.
$$

PROOF. For simplicity assume that the initial position $x \le x_1$. In this case the process $R_t^*$ as a solution to (4.2) is continuous. In view of (4.1) and (2.7),

$$
(4.4) \quad L^{a^*(R_s^*)}V(R_s^*) - M^*(R_s^*)V'(R_s^*) + M^*(R_t^*) = 0
$$

[since $M(1 - V'(x))^+ = M^*(x)(1 - V'(x))$], where the operator $L^a$ is defined in the proof of Theorem 2.2. Repeating the arguments of the proof of Theorem 2.2 and applying (4.4), we see that (we write $\tau$ instead of $\tau^\pi$ below, since there will be no confusion)

$$
(4.5) \quad E(e^{-\gamma(t \wedge \tau)}V(R_{t \wedge \tau}^*)) = V(x) - E\int_0^{t \wedge \tau} e^{-\gamma s}c_s^*\,ds.
$$

Taking the limit as $t \to \infty$, and applying (2.13), we obtain the desired result. $\square$

Combining Theorems 2.2 and 4.1, we get the following result immediately.

COROLLARY 4.2. *The function $V$ presented in Section 3 is the optimal return function and $\pi^*$ is the optimal policy.*

All the results we obtained are summarized in Tables 1 and 2 for easy reference.

Table 1
*The case of $\frac{2\delta}{\mu} < \alpha$*

| Range for $M$ | $x_\alpha$ | $x_\beta$ | $a^*(x)$ | Risk $\alpha$ ever attained | Risk $\beta$ ever attained | $x_1$ | $x_1$ is the first point at which the possible maximal risk is attained |
|---|---|---|---|---|---|---|---|
| $M > M_\beta$ | Positive and finite; see (3.28) | | (i) $\alpha$ for $x \in [0, x_\alpha)$ (ii) Increases from $\alpha$ to $\beta$ on $[x_\alpha, x_\beta]$; see (3.29) (iii) $\beta$ for $x \geq x_\beta$ | Yes | Yes | Positive; see (3.40) | No |
| $\frac{M_\alpha < M < M_\beta}{M = M_\beta}$ | Positive and finite; see (3.24) | $\infty$ | (i) $\alpha$ for $x \in [0, x_\alpha)$ (ii) Increases from $\alpha$ to $\frac{2\mu(\delta+M)}{\mu^2+2\gamma\sigma^2}$ on $[x_\alpha, x_1]$; see (3.29) (iii) $\frac{2\mu(\delta+M)}{\mu^2+2\gamma\sigma^2}$ for $x \geq x_1$ | Yes | $\frac{\text{No}}{\text{Yes}}$ | Positive; see (3.42) | Yes |
| $M_0(\alpha) < M \leq M_\alpha$ | $\infty$ | $\infty$ | $\alpha$ | Yes | No | Positive; see (3.52) | No |
| $M \leq M_0(\alpha)$ | $\infty$ | $\infty$ | $\alpha$ | Yes | No | 0 | Yes |

TABLE 2

| Range for $M$ | $x_\alpha$ | $x_\beta$ | $a^*(x)$ | Risk $\alpha$ ever attained | Risk $\beta$ ever attained | $x_1$ | $x_1$ is the first point at which the possible maximal risk is attained |
|---|---|---|---|---|---|---|---|
| | | | **The case of $\alpha \leq \frac{2\delta}{\mu} < \beta$** | | | | |
| $M > M_\beta$ | 0 | Positive and finite; see (3.55) | (i) Increases from $\frac{2\delta}{\mu}$ to $\beta$ on $[0,x_\beta]$; see (3.56); (ii) $\beta$ for $x > x_\beta$ | Yes if $\alpha = \frac{2\delta}{\mu}$; no if $\alpha > \frac{2\delta}{\mu}$ | Yes | Positive; see (3.40) | No |
| $\frac{M_0(2\delta/\mu)<M<M_\beta}{M=M_\beta}$ | 0 | $\infty$ | (i) Increases from $\frac{2\delta}{\mu}$ to $\frac{2\mu(\delta+M)}{\mu^2+2\gamma\sigma^2}$ on $[0,x_1]$; see (3.56) (ii) $\frac{2\mu(\delta+M)}{\mu^2+2\gamma\sigma^2}$ for $x \geq x_1$ | Yes if $\alpha = \frac{2\delta}{\mu}$; no if $\alpha > \frac{2\delta}{\mu}$ | $\frac{\text{No}}{\text{Yes}}$ | Positive; see (3.64) | Yes |
| $M_\alpha < M \leq M_0(\frac{2\delta}{\mu})$ | 0 | $\infty$ | $\frac{2\mu(\delta+M)}{\mu^2+2\gamma\sigma^2}$ | No | No | 0 | Yes |
| $M \leq M_\alpha$ | 0 | $\infty$ | $\alpha$ | Yes | No | 0 | Yes |
| | | | **The case of $\frac{\delta}{\mu} < \beta \leq \frac{2\delta}{\mu}$** | | | | |
| $M > M_0(\beta)$ | 0 | 0 | $\beta$ | No | Yes | Positive; see (3.72) | No |
| $M_\beta \leq M \leq M_0(\beta)$ | 0 | 0 | $\beta$ | No | Yes | 0 | Yes |
| $\frac{M_\alpha<M<M_\beta}{M=M_\alpha}$ | 0 | $\infty$ | $\frac{2\mu(\delta+M)}{\mu^2+2\gamma\sigma^2}$ | $\frac{\text{No}}{\text{Yes}}$ | No | 0 | Yes |
| $M < M_\alpha$ | $\infty$ | $\infty$ | $\alpha$ | No | No | 0 | Yes |

**5. Economic interpretation and conclusions.** The optimal policies obtained in the previous sections have clear economic meaning and are very easy to implement. Let us now elaborate.

The risk control policy is characterized by two critical reserve levels: $x_\alpha$ and $x_\beta$. The values of these two levels are further determined by four parameters: the minimum risk allowed ($\alpha$), the maximum risk allowed ($\beta$), the ratio between the debt rate and profit rate ($\frac{\delta}{\mu}$), and the maximum dividend rate allowed ($M$). Specifically, there are three different cases to consider.

The first case is when the company has very little debt compared to the potential profit (so that $\frac{2\delta}{\mu} < \alpha$). In this case, if the maximum dividend rate $M$ is large enough ($M > M_\beta$), then both critical reserve levels, $x_\alpha$ and $x_\beta$, are positive and finite. In other words, the company will minimize business activity (i.e., take the minimum risk $\alpha$) when the reserve is below $x_\alpha$, then gradually increase business activity when the reserve is between $x_\alpha$ and $x_\beta$, and then maximize business activity (i.e., take the maximum risk $\beta$) when the reserve reaches or goes beyond level $x_\beta$. This policy is the same as that obtained in [4] for the case of unbounded dividend rate. Next, if the maximum dividend rate $M$ is at a medium level ($M_\alpha < M < M_\beta$), then $x_\alpha$ remains positive and finite while $x_\beta$ becomes infinite. This implies that the company will become less aggressive; in particular, it will never take the maximum risk, due to a more restrictive dividend payout upper bound. Finally, if $M$ is so small that $M \le M_\alpha$, then both $x_\alpha$ and $x_\beta$ turn out to be infinite, meaning that business activities will be carried out at the minimum level or those business activities are redundant.

The second case is when the company has a higher debt–profit ratio (so that $\alpha \le \frac{2\delta}{\mu} < \beta$). In this case, $x_\alpha = 0$. This means that no matter how small the reserve is the company will never take the minimum risk; rather it will start with a bit higher risk level and gradually increase it. On the other hand, whether it will ever increase to the maximum possible risk (i.e., whether $x_\beta$ is finite or infinite) depends on the value of the maximum possible dividend rate $M$, in the same way as in the first case discussed above. Therefore, in the second case the company overall has to be a bit more aggressive than in the first case. This can be explained by the fact that when the debt rate is high, one needs to gamble on higher potential profits to get out of the "bankruptcy zone" as fast as possible, even at the expense of assuming higher risk.

The company becomes even more aggressive in the third case when the debt–profit ratio is even higher (precisely when $\frac{\delta}{\mu} < \beta \le \frac{2\delta}{\mu}$). In this case, when the maximum dividend rate $M$ is large enough ($M \ge M_\beta$), the maximum allowable risk $\beta$ is taken throughout, while the two critical levels $x_\alpha$ and $x_\beta$ are both zero. On the other hand, when $M$ is small enough so that $M < M_\alpha$, business activities are carried out at the minimum level $\alpha$ throughout.

On the other hand, the optimal dividend policy is always of a threshold type here the threshold is $x_1$ (which is positive or zero). Namely, the dividend distribution takes place only when the reserve exceeds the critical level $x_1$, in which case the dividend payout rate is $M$.

It is interesting to note that in the case of unbounded dividend rate, the maximum business activity is always taken *before* dividend distributions ever take place; see [4]. However, in the present case of bounded dividend rate, the company may need to pay dividends *before* the maximum risk level $\beta$ is ever taken; refer to Tables 1 and 2 for details. This represents a striking difference between the cases of unbounded and bounded dividend rate. The economic reason for such behavior is the following. When there is a significant constraint on the dividend rate, there may be no necessity to pursue business aggressively because the accumulated liquid assets cannot be paid out as dividends fast enough anyway.

In conclusion, we point out an intricate interplay between the restriction on the dividend distribution rate and that on the risk control of a financial company. The sheer number of qualitatively different optimal policies, which appears due to different possible relationships between exogenous parameters, shows the multiplicity of different economic environments which a financial company faces depending on the size of the debt, on the constraint on the dividend rate, and on the size of available business activity.

T. CHOULLI
MATHEMATICAL AND
    STATISTICAL SCIENCES DEPARTMENT
UNIVERSITY OF ALBERTA
EDMONTON, ALBERTA
CANADA T6G 2G1

M. TAKSAR
DEPARTMENT OF MATHEMATICS
UNIVERSITY OF MISSOURI
COLUMBIA, MISSOURI 65211
USA
E-MAIL: taksar@math.missouri.edu

X. Y. ZHOU
DEPARTMENT OF SYSTEMS ENGINEERING
    AND ENGINEERING MANAGEMENT
CHINESE UNIVERSITY OF HONG KONG
SHATIN, HONG KONG
PEOPLES REPUBLIC OF CHINA
E-MAIL: xyzhou@se.cuhk.edu.hk